\documentclass[12pt]{article}
\usepackage{authblk}
\usepackage{abstract}
\usepackage{cite}
\usepackage{figsize,subfigure,graphicx,ifthen,calc,psfrag,epsfig}
\usepackage{longtable}
\usepackage{multirow}
\usepackage{ulem}
\usepackage{rotating}
\usepackage{graphics}
\usepackage{graphicx}
\usepackage{epstopdf}
\usepackage{epsfig}
\usepackage{indentfirst}
\usepackage{amssymb,amsmath}
\usepackage{amsmath}
\usepackage{color}
\usepackage{amsfonts}
\usepackage{amssymb}
\usepackage{subfigure}
\usepackage{latexsym,bm}
\usepackage{graphicx}
\usepackage{amsmath}
\usepackage{algorithm,algorithmic}
\usepackage{diagbox}
\usepackage{appendix}
\usepackage{listings}
\newtheorem{example}{Example}

\begin{document}
\bibliographystyle{elsarticle-num}

\title{Physics-informed Neural Networks for Elliptic Partial Differential Equations on 3D Manifolds}

\author[a,b]{Zhuochao Tang}
\author[a,b,c]{Zhuojia Fu \thanks{Corresponding author, email address: paul212063@hhu.edu.cn (Zhuojia Fu).}}
\affil[a]{Key Laboratory of Ministry of Education for Coastal Disaster and Protection, Hohai University, Nanjing 210098, China}
\affil[b]{Center for Numerical Simulation Software in Engineering and Sciences, College of Mechanics and Materials, Hohai University, Nanjing 211100, China}
\affil[c]{State Key Laboratory of Mechanics and Control of Mechanical Structures, Nanjing University of Aeronautics and Astronautics, Nanjing 210098, China}
\renewcommand*{\Affilfont}{\small\it}
\renewcommand\Authands{ and }
\date{}
\maketitle

\begin{abstract}
Motivated by recent research on Physics-Informed Neural Networks (PINNs), we make the first attempt to introduce the PINNs for numerical simulation of the elliptic Partial Differential Equations (PDEs) on 3D manifolds. PINNs are one of the deep learning-based techniques. Based on the data and physical models, PINNs introduce the standard feedforward neural networks (NNs) to approximate the solutions to the PDE systems. By using automatic differentiation, the PDEs system could be explicitly encoded into NNs and consequently, the sum of mean squared residuals from PDEs could be minimized with respect to the NN parameters. In this study, the residual in the loss function could be constructed validly by using the automatic differentiation because of the relationship between the surface differential operators $\nabla_S/\Delta_S$ and the standard Euclidean differential operators $\nabla/\Delta$. We first consider the unit sphere as surface to investigate the numerical accuracy and convergence of the PINNs with different training example sizes and the depth of the NNs. Another examples are provided with different complex manifolds to verify the robustness of the PINNs.
\end{abstract}


\section{Introduction}
Various applications in science and engineering, as a matter of fact, refer to numerical solutions of Partial Differential Equations (PDEs) on surfaces or more general manifolds. Such applications include the generation of textures\cite{Witkin1995Reaction} or the visualization of vector fields\cite{2000Anisotropic} in image processing, flows and solidification\cite{2004A} on surfaces in fluid dynamics and evolving surfactants\cite{2003An} on interfaces in biology, etc. Typical numerical methods solve the surface PDEs by means of either intrinsic\cite{2005Surface}, extrinsic\cite{Marcelo2001Variational}, or embedding techniques\cite{2008A}. Both intrinsic and extrinsic methods require analytical expressions of differential operators defined on manifolds, which are mapped from standard ones defined in Cartesian space. The difference between intrinsic and extrinsic methods is that intrinsic methods need parameterization on the surface while extrinsic methods impose nodes/meshes on manifolds. On the other side, the embedding techniques extend the PDE on manifolds into one additional embedding space and then solve it. In recent years, various methods including either mesh-based methods like the Finite Element Method (FEM)\cite{2014A,2014Error}, Finite Volume Method (FVM)\cite{2009A} or mesh-free methods like RBF-type methods\cite{2019Kernel,2019Kernel1,2013A}, meshless finite difference methods\cite{Pratik2019A} have been proposed to solve PDEs on manifolds. Different from these typical numerical methods, machine learning methods\cite{Goldberg1989Genetic} provide another option for solving PDEs on manifolds. 

Machine learning approach could be dated back to last century while it has not been widely taken into practice in engineering applications until recent years, for instance, computer vision\cite{2006Tensor}, natural language processing\cite{2017Integrating} and face recognition\cite{2015Cloud}, etc. Such approaches\cite{2015Deep} using deep neural networks (NNs) build a black-box model based on "training data" to make predictions and recently, they have also succeeded in solving Partial Differential Equations (PDEs) system. Physics-Informed Neural Networks (PINNs)\cite{2018Physics,2020Conservative}, as one of these deep learning methods, are firstly put forward in \cite{0Hidden} and have been applied to many different problems such as fractional PDEs\cite{2019fPINNs}, computational fluid dynamics\cite{2020Physics} and multi-physics areas\cite{2020Physics1}, etc. PINNs aim to replace the PDE solution with a neural network and take advantage of information from PDEs and initial/boundary conditions to form an optimized system explicitly. This explicit system origins from the information based on training data and could also be defined as terminology loss function. By minimizing this system with respect to the parameters defined in NNs, PINNs could find one NN which best describes the physical behaviour governed by the PDEs. To specify the differential operators acting on the variables, PINNs employ the automatic differentiation technology and classical chain rule. Motivated by the recent development of the PINNs, we make the first attempt to utilize this novel method, Physics-Informed Neural Networks (PINNs), to solve the second-order elliptic Partial Differential Equations (PDEs) on 3D manifolds.

Benefit from the analytical expressions between surface differential operators and Euclidean differential operators in intrinsic/extrinsic methods, PINNs could be connected seamlessly with it to solve PDEs on manifolds. Another advantage referred here is that PINNs nearly have no restrictions on the distribution of data points. For classical methods, the higher the quality of the mesh or point cloud, the more computing resources are consumed while PINNs avoid this. This also shows a good flexibility and potential of PINNs to deal with the high-dimensional problems with complex shapes.

The remaining paper is organized as follows: Section 2 gives details on PDEs defined on manifolds, introduces the PINNs and describes its implementation. In section 3, we demonstrate the effectiveness of PINNs under several numerical examples. In this section, we first illustrate the convergence results by using different parameters in PINNs and test the robustness of PINNs by adopting sundry smooth manifolds. In the same section, we also make a comparison of numerical results by using arbitrary training points and quasi-uniform points. Some conclusions are summarized in Section 4.

\section{Methodology}
In this section, we give a detailed introduction of surface operators and explanation of how PINNs solve second-order elliptic partial differential equations on 3D manifolds.

\subsection{Partial Differential Equations (PDEs) on 3D manifolds}
In this study, we pay our attention to second-order elliptic partial differential equations (PDEs) posed on some sufficiently smooth, connected, and compact surface $S\subset\mathbb{R}$ with no boundary and $dim(S)=d-1$. We will focus on the case of $d=3$ for notational simplicity in the following description. Any other cases with different $d$ could be extended obviously. We first consider general second-order elliptic partial differential equations
\begin{equation}
 \label{eq1}
\left(a\Delta_S-\vec{\mathbf{b}}\cdot\nabla_S+c\right)u(x,y,z)=f(x,y,z) 
\end{equation}
with the certain coefficients $a, \vec{\mathbf{b}}, c$ and $\Delta_S, \nabla_S$ are respectively the Laplace-Beltrami (a.k.a. surface Laplacian) operator and the surface gradient operator.

The main difference between surface PDEs defined on manifolds and standard PDEs posed in some bounded domains with flat geometry is that the curvatures of manifold play vital roles in physical behaviours governed by the PDEs. To describe the relationship between surface operators and standard operators, we denote the unit outward normal vector at any $\mathbf{x}\in S$ as $\mathbf{n}=(n^x, n^y, n^z)$ and the corresponding projection matrix to the tangent space as
\begin{equation}
\mathbf{P}(x)=(\mathbf{I}_3-\mathbf{nn}^T)\in \mathbb{R}^{3\times 3}, \label{eq2}
\end{equation}
where $\mathbf{I}$ is the 3-by-3 identity matrix. Then the surface gradient operator $\nabla_S$ could be defined in terms of the standard Euclidean gradient $\nabla$ via projections as
\begin{equation}
\nabla_S:=\mathbf{P}\nabla, \label{eq3}
\end{equation}
and similarly, the Laplace-Beltrami operator $\Delta_S$ could be defined as
\begin{equation}
\Delta_S:=\nabla_S\cdot\nabla_S. \label{eq4}
\end{equation}

As mentioned in Section 1, direct approximation of surface differential operators is not an easy thing because the analytic expression of surface Laplace-Beltrami operator involves the derivatives of normal vector. In \cite{2018AA}, one significant formula between surface differential operators and Euclidean differential operators acting on any sufficiently smooth function could be found:
\begin{equation}
\nabla_Su:=\nabla_Su-\mathbf{n}\partial_{\mathbf{n}}u, \label{eq5}
\end{equation}
\begin{equation}
\Delta_Su:=\Delta_Su-H_S\partial_{\mathbf{n}}u-\partial_{\mathbf{n}}^{(2)}u. \label{eq6}
\end{equation}
in which $\partial_{\mathbf{n}}u-\mathbf{n}^T\nabla u, \partial_{\mathbf{n}}^{(2)}u:=\mathbf{n}^TJ(\nabla u)\mathbf{n}$ and $H_S=trace\left(J(\mathbf{n}(\mathbf{I}-\mathbf{nn}^T))\right)$. Here, $J$ means the Jacobian operator in Euclidean space.

For better understanding, we give one example to derive the explicit expression of surface differential operators on the unit sphere. Simplifying with the manifold $S=x^2+y^2+z^2-1$, one could naturally obtain the unit normal vector $[x \quad y \quad z]^T$. Putting this into Eqs.\eqref{eq5}\eqref{eq6}, the surface differential operators in Cartesian coordinates are represented by
\begin{equation}\label{eq7}
\nabla_S=\left[
\begin{matrix}
1-x^2 & -xy   & -xz \\
-xy   & 1-y^2 & -yz \\
-xz   & -yz   & 1-z^2
\end{matrix}\right] \left[
\begin{matrix}
\partial_x \\ \partial_xy \\\partial_z
\end{matrix}\right] = \left[
\begin{matrix}
(1-x^2)\partial_x-xy\partial_y-xz\partial_z \\ -xy\partial_x+(1-y^2)\partial_y-yz\partial_z \\ -xz\partial_x-yz\partial_y+(1-z^2)\partial_z
\end{matrix}\right],
\end{equation}
\begin{equation}\label{eq8}
\Delta_S=(1-x^2)\partial_{xx}+(1-y^2)\partial_{yy}+(1-z^2)\partial_{zz}-2xy\partial_{xy}-2xz\partial_{xz}-2yz\partial_{yz}-2x\partial_{x}-2y\partial_{y}-2z\partial_{z}.
\end{equation}
Once obtaining the Eqs.\eqref{eq7}\eqref{eq8}, any continuous second-order elliptic differential operators $\left(a\Delta_S-\vec{\mathbf{b}}\cdot\nabla_S+c\right)$ defined on smooth manifolds could be expressed by continuous Euclidean operators.

\subsection{Physics-Informed Neural Networks (PINNs)}
The main aim of PINNs is to approximate the solutions to PDEs. In this section, we focus on introducing the basic idea of PINNs and how it solves PDEs on manifolds. We consider the Eq.\eqref{eq1} again here:
\begin{equation}\label{eq9}
\left(a\Delta_S-\vec{\mathbf{b}}\cdot\nabla_S+c\right)u(x,y,z)=f(x,y,z).
\end{equation}

\begin{figure}[!ht]
\centering
\includegraphics[width=6in]{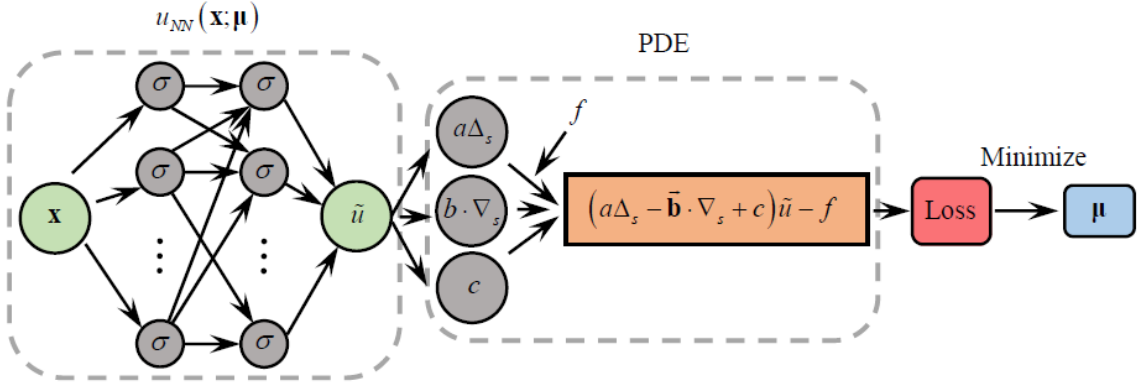}
\caption{Sketch of Physics-Informed Neural Networks (PINNs): basic process of PINNs for PDEs on manifolds.} \label{fig1}
\end{figure}

In the PINNs, there are three different ways to construct the approximate solutions $u(x,y,z)$ to the PDEs\cite{2019fPINNs}. Due to fact that the PDE \eqref{eq9} defined on closed manifolds have no boundary conditions, the direct construction of the approximate solutions is employed in this work as an output of Neural Networks (NN), namely, $\tilde{u}(x,y,z)=u_{NN}(\mathbf{x};\bm{\mu}),\mathbf{x}\in S$. The NN, which is parameterized with finitely many weights and biases, acts as a substitute model of the PDE model to approximate the mapping from the spatial coordinates to the solutions of equation. As shown in Fig.\ref{fig1}, one NN usually contains multiple hidden layers to obtain more accurate solutions. Here, PINNs seek to optimize the NN parameters composed of weights $w$ and biases $b$ by minimizing so-called loss function. Usually, the loss function is defined as the sum of mean-squared-error from both governing equations (PDEs) and boundary conditions on the training points. For PDEs defined on manifolds without boundary conditions, the loss function is expressed as
\begin{equation} \label{eq10}
Loss(\bm{\mu})=\frac{1}{N}\sum_{k=1}^N\left[\left(a\Delta_S-\vec{\mathbf{b}}\cdot\nabla_S+c\right)\tilde{u}(\mathbf{x}_k)-f(\mathbf{x}_k)\right]^2,
\end{equation}
in which $N$ is the total number of the training points. By introducing the relationship between surface differential operator and Euclidean differential operator, the loss function finally could be parameterized by the NN parameter $\bm{\mu}$:
\begin{equation} \label{eq11}
Loss(\bm{\mu})=\frac{1}{N}\sum_{k=1}^N\left[\left(a(\Delta-H_S\partial_{\mathbf{n}}-\partial_{\mathbf{n}}^{(2)})-\vec{\mathbf{b}}\cdot(\nabla-\mathbf{n}\partial_{\mathbf{n}})+c\right)u_{NN}(\mathbf{x}_k;\bm{\mu})-f(\mathbf{x}_k)\right]^2.
\end{equation}

We minimize the loss function with respect to $\bm{\mu}$ and then we could obtain a minimizer. The automatic differentiation technique and the chain rule are used in loss function to compute the spatial derivatives of $u_{NN}(\mathbf{x}_k;\bm{\mu})$. This optimization process is called “training”. We have to emphasize here that for a moderate size of NN, there always exist infinitely many global minimizers and we can’t guarantee that the minimizer found coincides with the exact solution. However, for linear second-order elliptic and parabolic PDEs posed in some bounded domains, it has been proved that the sequence of minimizers strongly converges to the PDE solution in $C^0$ \cite{Shin_2020}, which still applies to linear second-order elliptic PDEs defined on manifolds. Nevertheless, we will still use multiple sets of initial NN parameters $\bm{\mu}$ in the following numerical examples to avoid its uncertainty.
To specify the superiority of PINNs, we would employ arbitrary irregular training points on 3D manifolds in our work and make specific comparisons.

\section{Numerical examples}
In this section, several different examples are provided. We first explore the convergence and the accuracy of PINNs for second-order elliptic PDEs on the unit sphere and then more 3D manifolds are taken into consideration to verify its robustness. To quantify the accuracy and effectiveness of our approximate solutions, we introduce the $L_2$ error measures as follows:
 \begin{equation}\label{eq12}
L_2 \ \text{error}= \frac{\sqrt{\sum\limits_{k=1}^{N}\left(u(\mathbf{x}_k)-\tilde{u}(\mathbf{x}_k)\right)^2}}{\sqrt{\sum\limits_{k=1}^{N}\left(u(\mathbf{x}_k)\right)^2}},
 \end{equation}
where $u(\mathbf{x}_k), \tilde{u}(\mathbf{x}_k)$ represent the reference solution and approximate solution at the $k$-th point. To avoid the uncertainty of different initializations for the network parameters $\bm{\mu}$ and find an optimal neural network as much as possible, we employ optimization method of ``L-BFGS'' and plot the mean for the solution errors from the 10 runs, whom we adopt as a new metric of convergence. Xavier initialization is taken into consideration and all the tests are implemented in Python by using Tensorflow.

\begin{example}
Convergence and accuracy test on unit sphere

\rm{In this example, the parameters in Eq.\eqref{eq1} are chosen as $a=1, \vec{\mathbf{b}}=[1 \ 1 \ 1]^T, c=5$ and the reference solution is assembled by trigonometric function, which is expressed as
\begin{equation} \label{eq13}
u(x,y,z)=\sin x\sin y\sin z.
\end{equation}
The force term is simply obtained by substituting reference solution into the equation. The points with total number of 2500 are set to be distributed on the unit sphere as shown in Fig.\ref{fig2}. Here, we select $N$ points randomly from these quasi-uniform points as training data and all these 2500 points are regarded as test points to test the convergence of PINNs. As derived above in Eqs.\eqref{eq7}\eqref{eq8}, the loss function on this unit sphere could be obtained easily. 
\begin{figure}[!ht]
\centering
\includegraphics[width=3in]{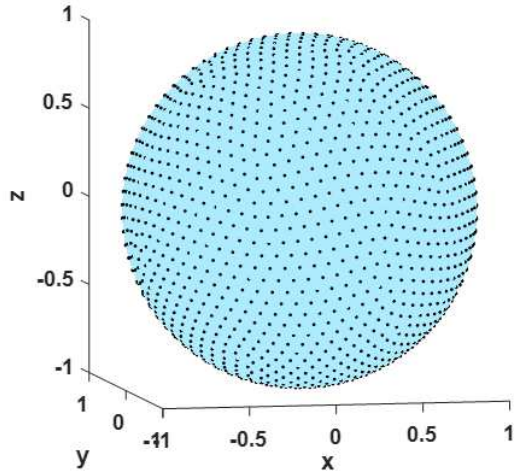}
\caption{Sketch of the quasi-uniform points distributed on the unit sphere: the point sets could be obtained by using minimum energy (ME) algorithm\cite{0Numerical}.} \label{fig2}
\end{figure}

Since we have no idea on how sensitive PINNs approximations are to surface differentiation operators, we attempt to use various NNs with different number of hidden layers (a.k.a. the depth of the NN, say, 4 hidden layers mean depth is 5) and neurons (a.k.a. the width of the NN, say, 20 neurons mean width is 20). Fig.\ref{fig3} shows the convergence results and Fig.\ref{fig4} indicates some snapshots of error distribution by using different parameters.
\begin{figure}[!ht]
\centering
\subfigure[]{\label{fig3a}
\includegraphics[scale=0.4]{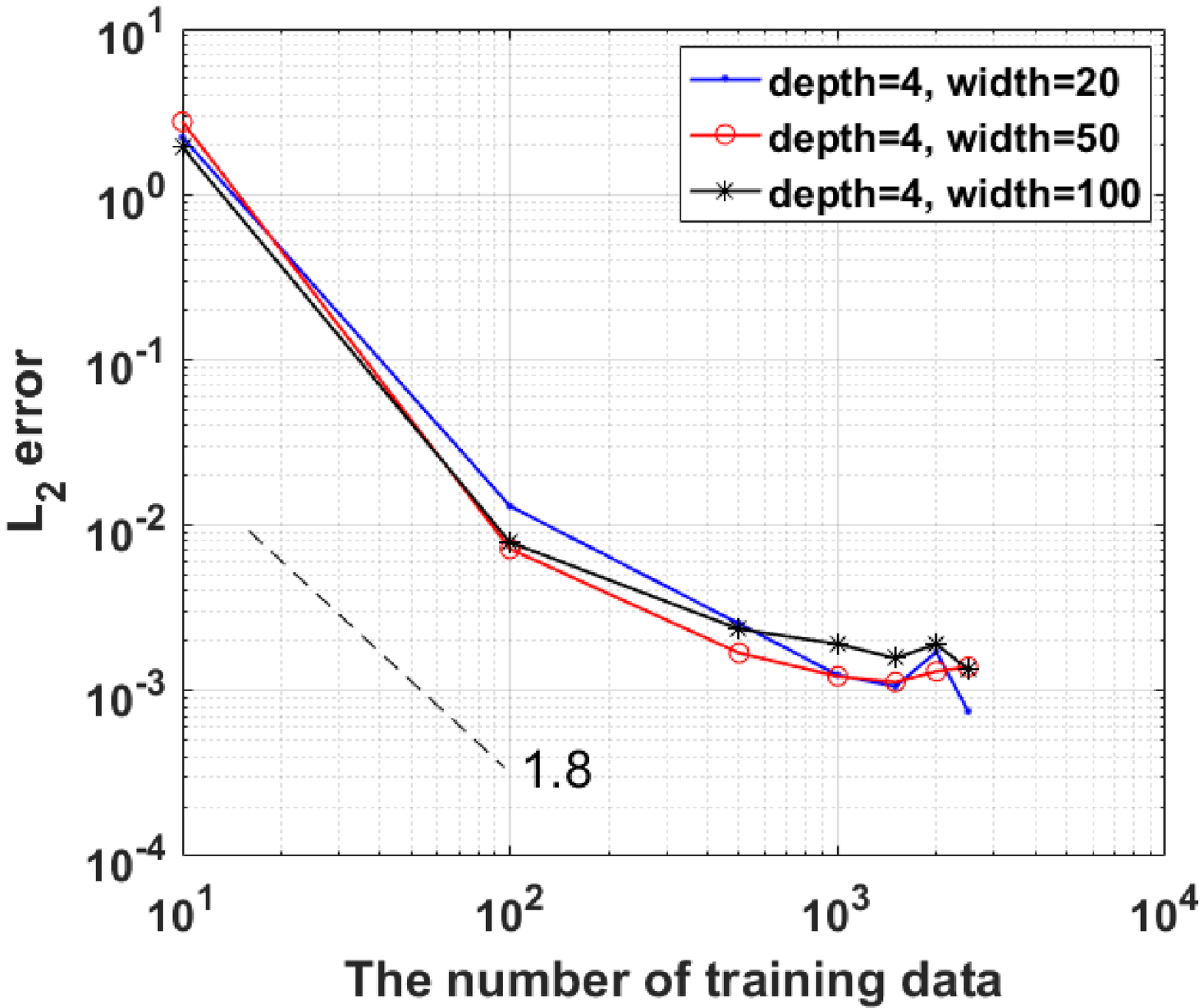}}
\subfigure[]{\label{fig3b}
\includegraphics[scale=0.4]{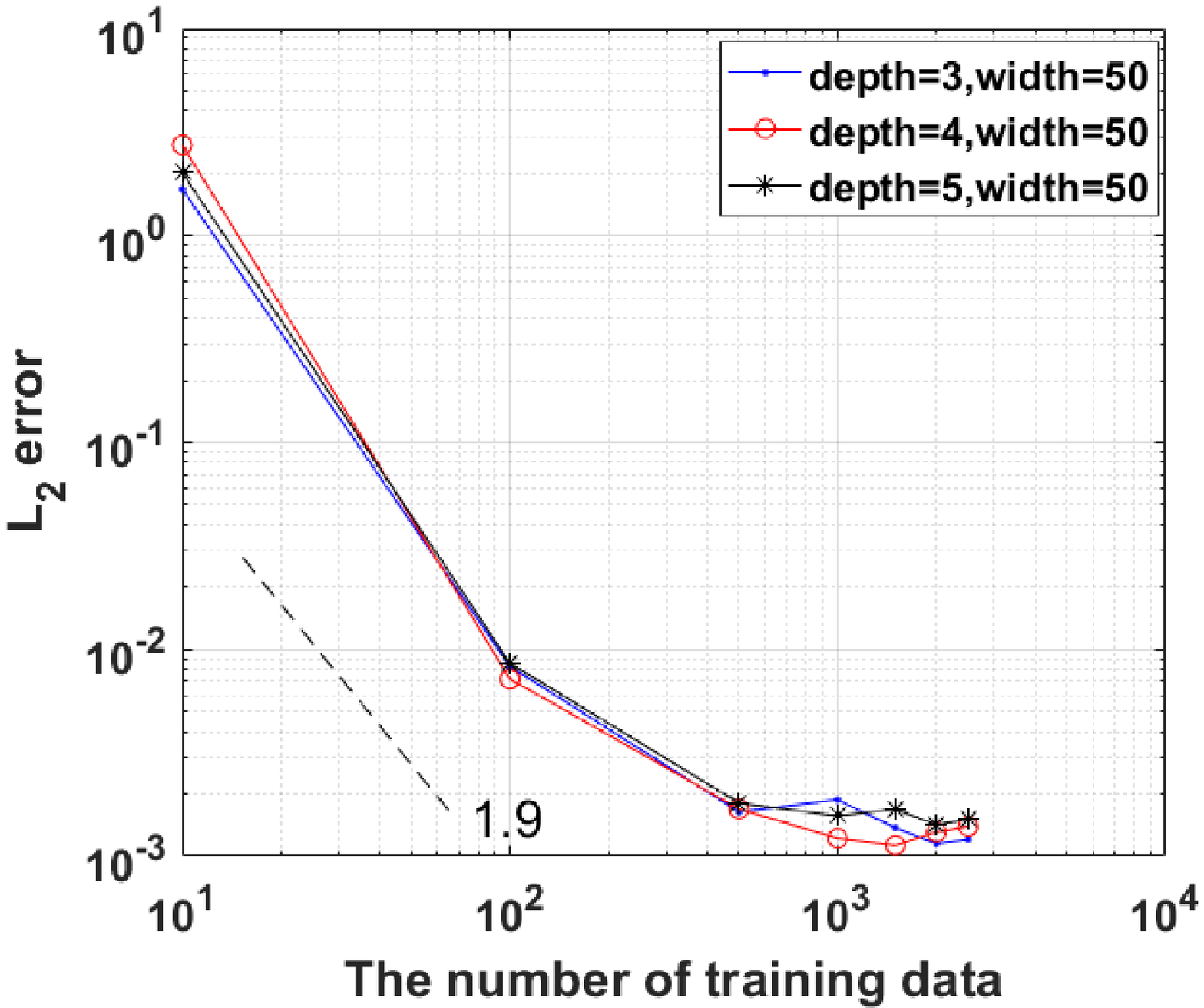}}
\caption{Example 1: Convergence results by using (a) different widths and (b) different depths.} \label{fig3}
\end{figure}
\begin{figure}[!ht]
\centering
\subfigure[]{\label{Fig4a}
\includegraphics[scale=0.4]{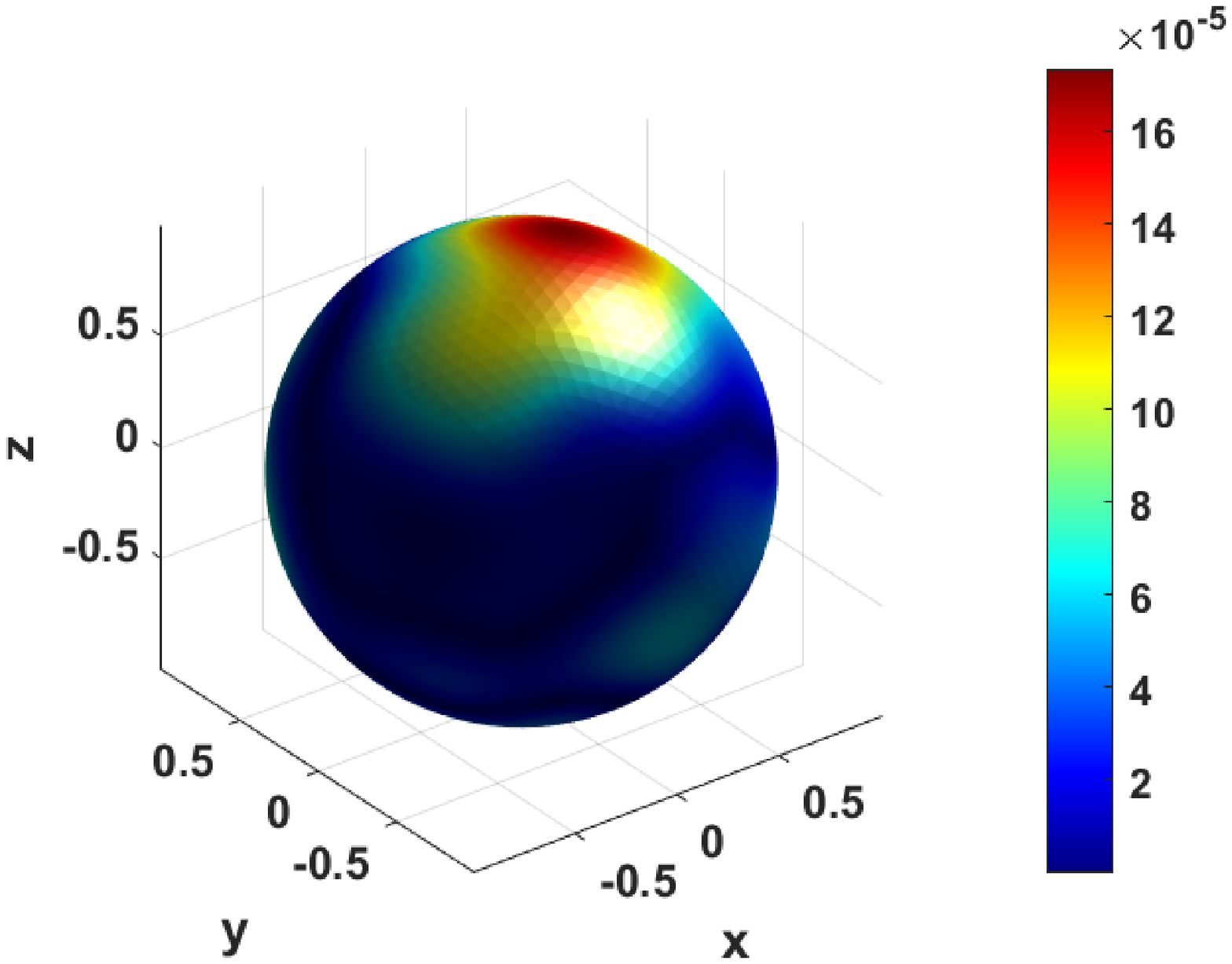}}
\subfigure[]{\label{Fig4b}
\includegraphics[scale=0.4]{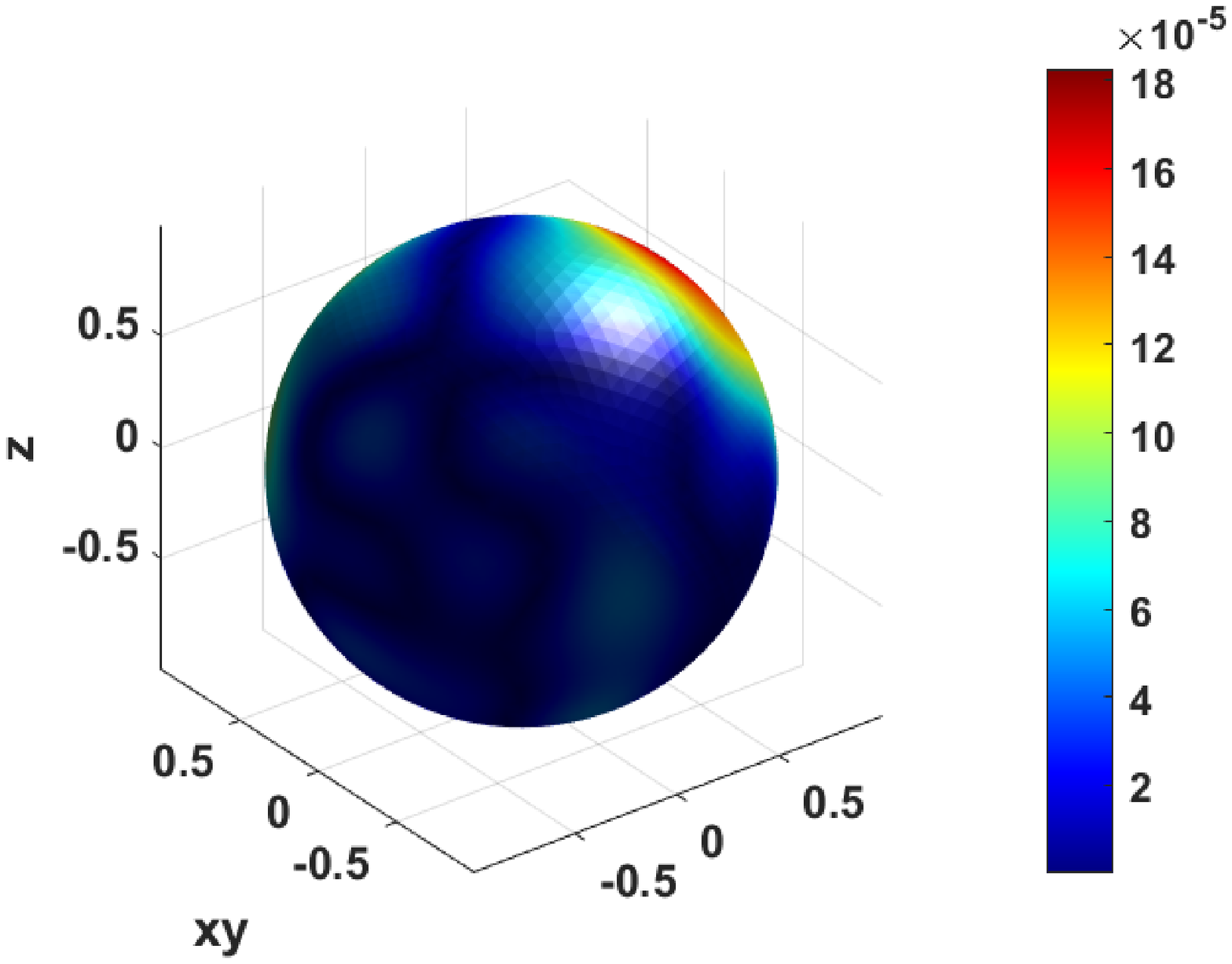}}
\subfigure[]{\label{Fig4c}
\includegraphics[scale=0.4]{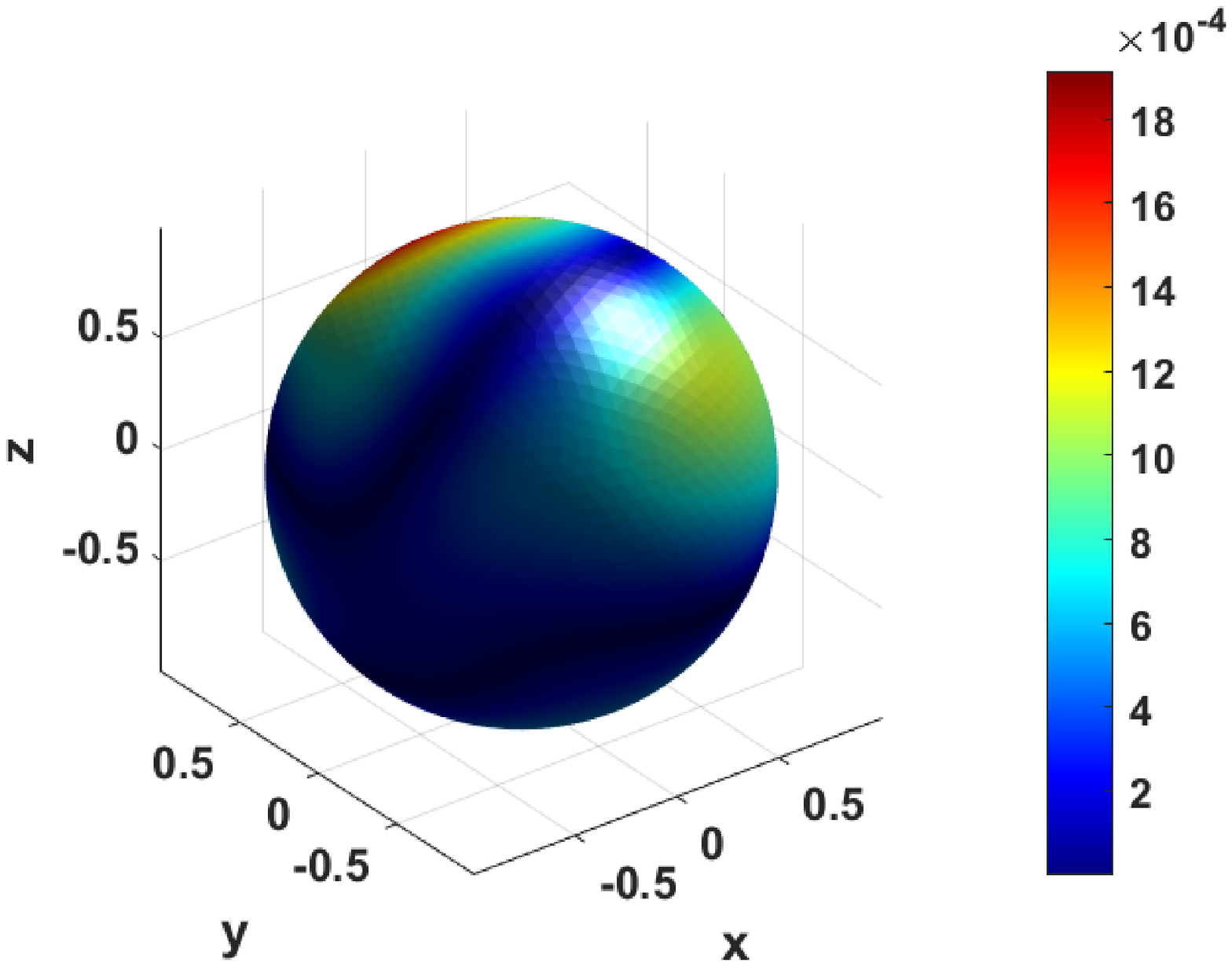}}
\subfigure[]{\label{Fig4d}
\includegraphics[scale=0.4]{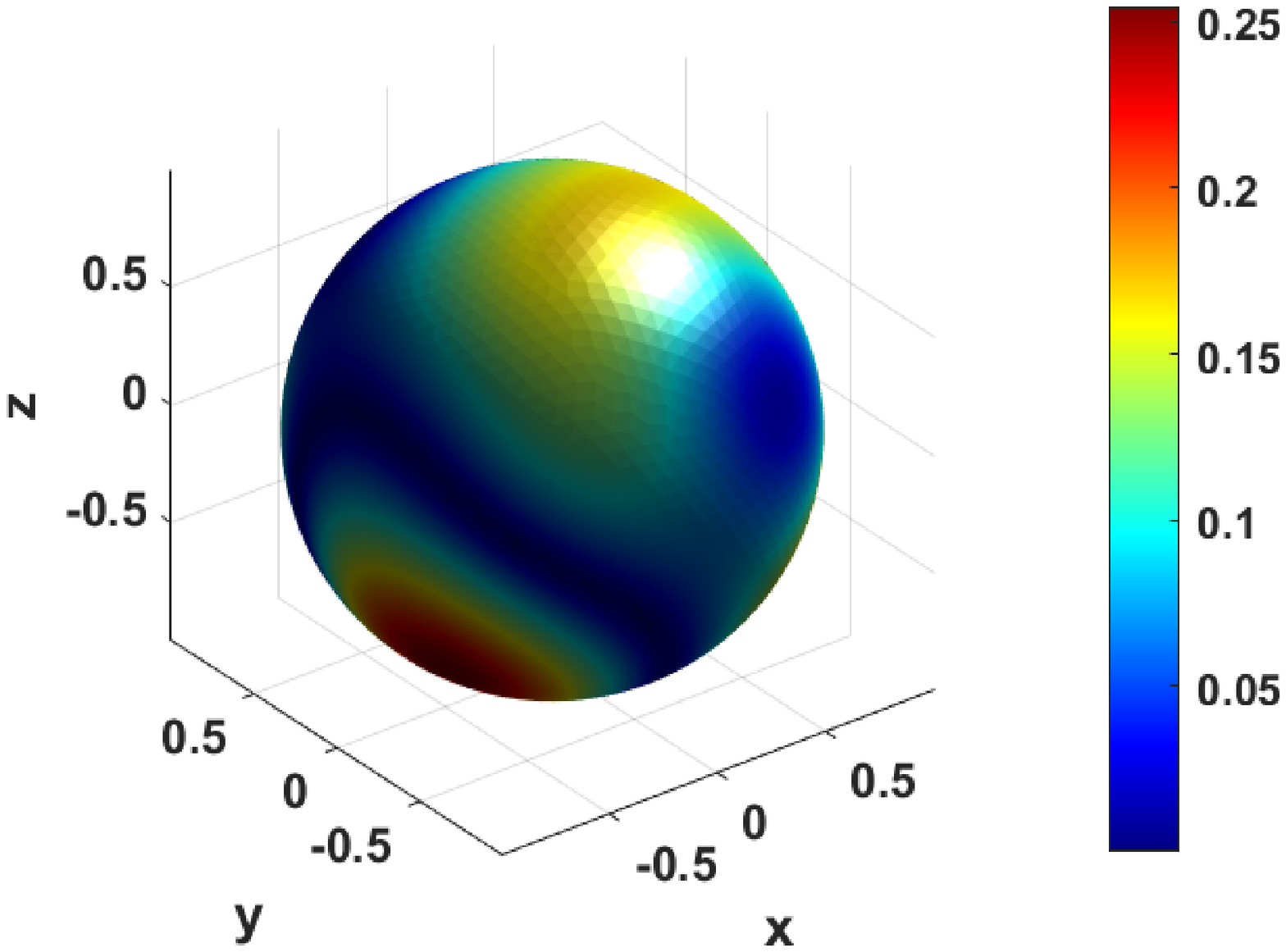}}
\caption{Example 1: Random few snapshots of error distribution with 3 hidden layers and 50 neurons by using (a) 2500 training data (b) 1500 training data (c) 100 training data (d) 10 training data.}\label{fig4}
\end{figure}

As seen from Figs.\ref{fig3} and \ref{fig4}, we use respectively 2, 3, 4 hidden layers with 20, 50, 100 neurons to test convergence and accuracy of PINNs in solving Eq.\eqref{eq1}. Numerical results converge to around $10^{-4}\sim 10^{-3}$ with convergence rate of 1.8 and 1.9. Apparently, we could obtain similar results form different parameters. This shows PINNs approximation have good adaptability to surface differential operators. In order to confront the following more complex manifolds, we set 3 hidden layers with 50 neurons by default in PINNs. Furthermore, we particularly plot the error distribution in Fig.\ref{fig4} to visualize the results, which shows a good accuracy of PINNs explicitly in solving surface PDEs.
}
\end{example}

\begin{example}
PDEs on complex manifolds

\rm{In this example, we attempt to test the robustness of PINNs in solving PDEs on more complex manifolds and make a direct comparison by using quasi-uniform training data and by using randomly distributed training data. The parametric equations or implicit expressions of some manifolds used in this example including Torus, a constant distance product (CDP) manifold, Bretzel2, Orthocircle, Red Blood Cell (RBC) are provided here:
\begin{itemize}
  \item [1)] Tours: $S=\left(1-\sqrt{x^2+y^2}\right)^2+z^2-\frac{1}{9}$;       
  \item [2)] CDP:
  $S=\sqrt{(x-1)^2+y^2+z^2}\sqrt{(x+1)^2+y^2+z^2}\sqrt{x^2+(y-1)^2+z^2}\\
\qquad \sqrt{x^2+(y+1)^2+z^2}-1.1$;
  \item [3)] Bretzel2: $S=\left(x^2(1-x^2)-y^2\right)^2+\frac{1}{2}z^2-\frac{1}{40}$;
  \item [4)] Orthocircle: $S=\left((x^2+y^2-1)^2+z^2\right)\left((y^2+z^2-1)^2+x^2\right)\left((z^2+x^2-1)^2+y^2\right)- \\ 0.075^2\left(1+3(x^2+y^2+z^2)\right)$;
  \item [5)] RBC: 
\begin{equation} \nonumber
S=\left\{
\begin{array}{lr}
      x=1.15\cos(\lambda)\cos(\theta), &  \\
      y=1.15\sin(\lambda)\cos(\theta), & -\pi\leq \lambda \leq \pi, \frac{-\pi}{2}\leq \theta \leq \frac{\pi}{2}.\\
      z=0.5\sin(\lambda)\left(0.24+2.3\cos(\theta)^2-1.3\cos(\theta)^4\right), &  
\end{array}
\right.
\end{equation}
\end{itemize}

In this test, the coefficients in Eq.\ref{eq1} are set as $a=1, \vec{\mathbf{b}}=[1 \ 1 \ 1]^T, c=1$ and the reference solution is changed to $u(x,y,z)=\sin x\cos y \sin z$. We first employ Torus manifold by using 500 quasi-uniform training data and by using 500 randomly distributed training data to make a comparison.
\begin{figure}[!ht]
\centering
\subfigure[]{\label{fig5a}
\includegraphics[scale=0.4]{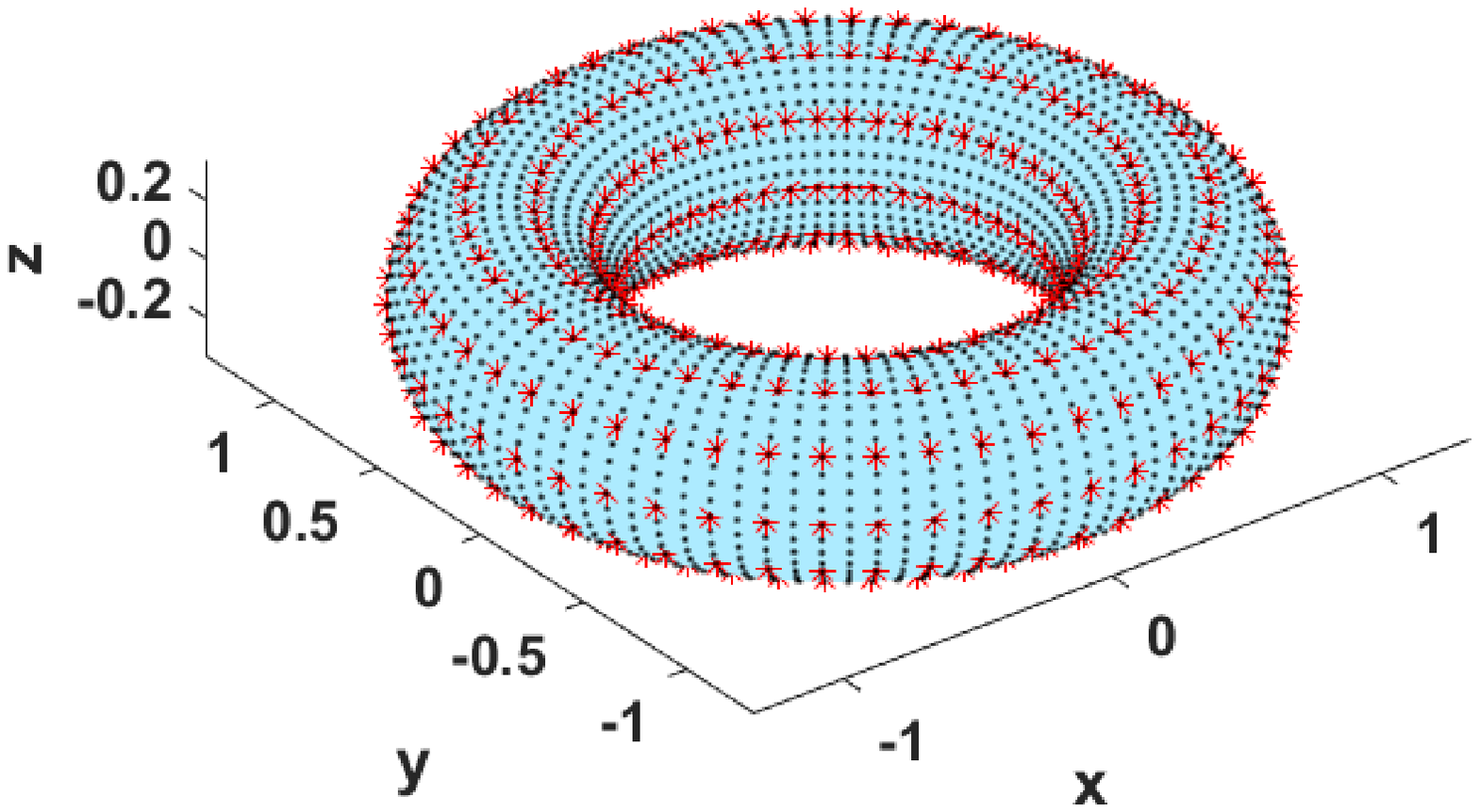}}
\subfigure[]{\label{fig5b}
\includegraphics[scale=0.4]{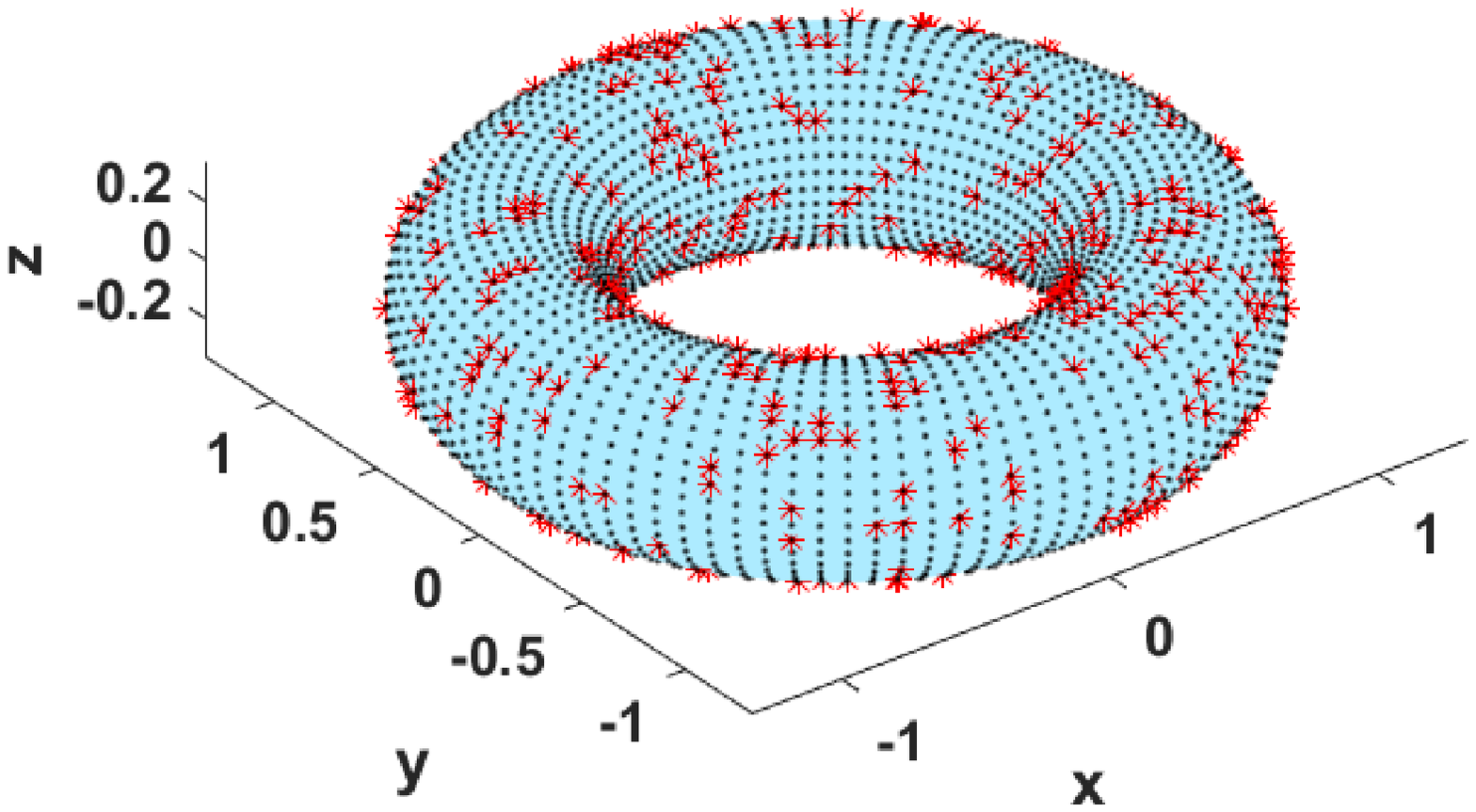}}
\caption{Example 2: Two different selections of training data on Torus manifold (a) quasi-uniform training data (b) random training data: red ``*'' points are selected training data; black points are test points.} \label{fig5}
\end{figure}
\begin{figure}[!ht]
\centering
\subfigure[]{\label{fig6a}
\includegraphics[scale=0.4]{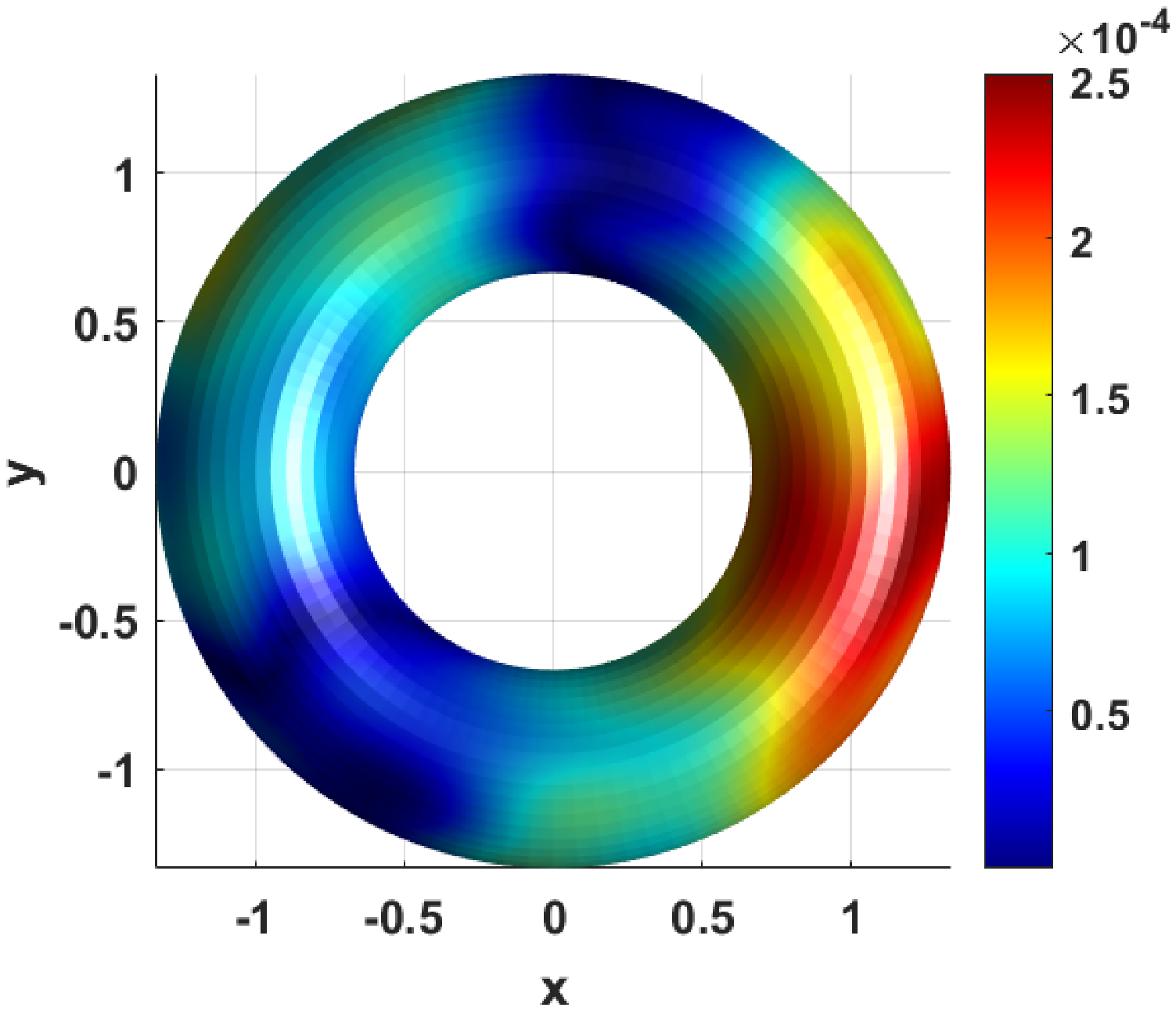}}
\subfigure[]{\label{fig6b}
\includegraphics[scale=0.4]{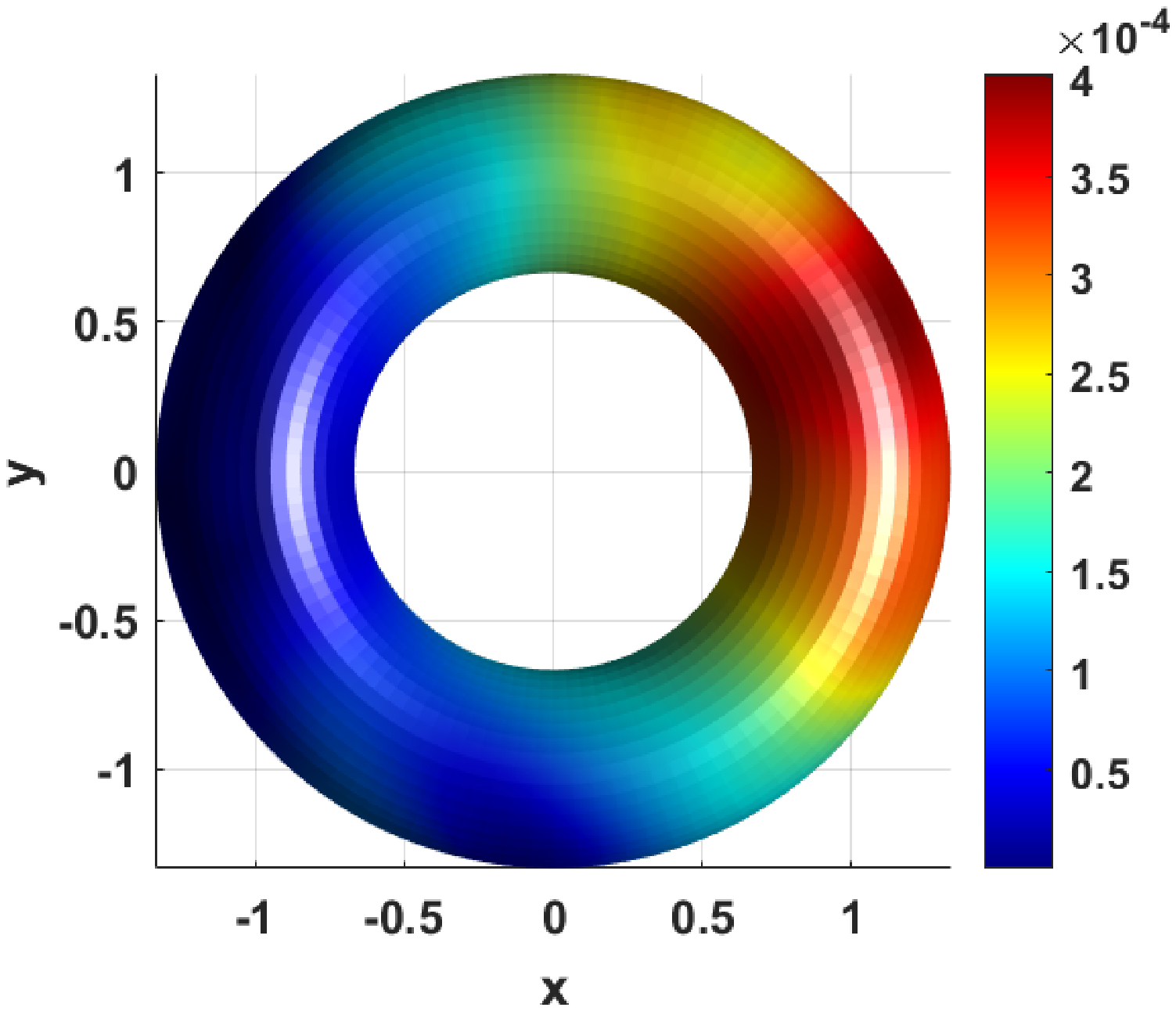}}
\caption{Example 2: Snapshots of error distribution (a) by using quasi-uniform training data and (b) by using randomly distributed training data.} \label{fig6}
\end{figure}

It can be found from Fig.\ref{fig6} that the distribution of training data hardly affect the numerical results, which proves that PINNs nearly have no restrictions on distribution of training points as long as the data themself are accurate enough. This shows again that PINNs to some extent are superior to some typical numerical methods in high dimensional problems. In addition, distribution numerical errors and $L_2$ errors on different manifolds are given respectively in Fig.\ref{fig7} and Tab.\ref{tab1}. The number of training points is chosen as 500 and the total number of points corresponding to CDP, Breztel2, Orthocircle and RBC is 3996, 3690, 4286 and 4000. Dealing with high-dimensional PDEs defined on complex manifolds, PINNs show a good stability and robustness.
\begin{figure}[!ht]
\centering
\subfigure[]{\label{Fig7a}
\includegraphics[scale=0.4]{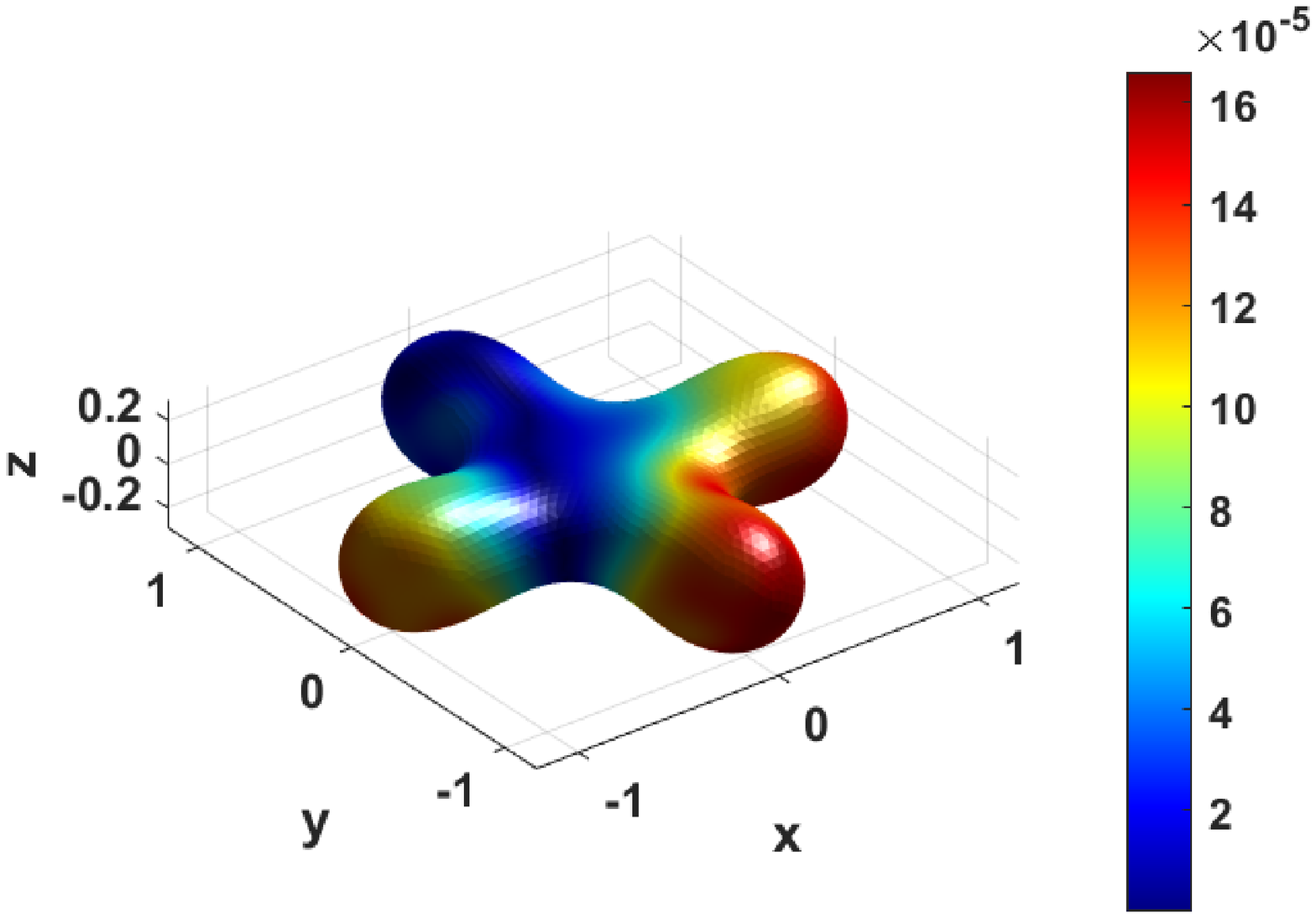}}
\subfigure[]{\label{Fig7b}
\includegraphics[scale=0.4]{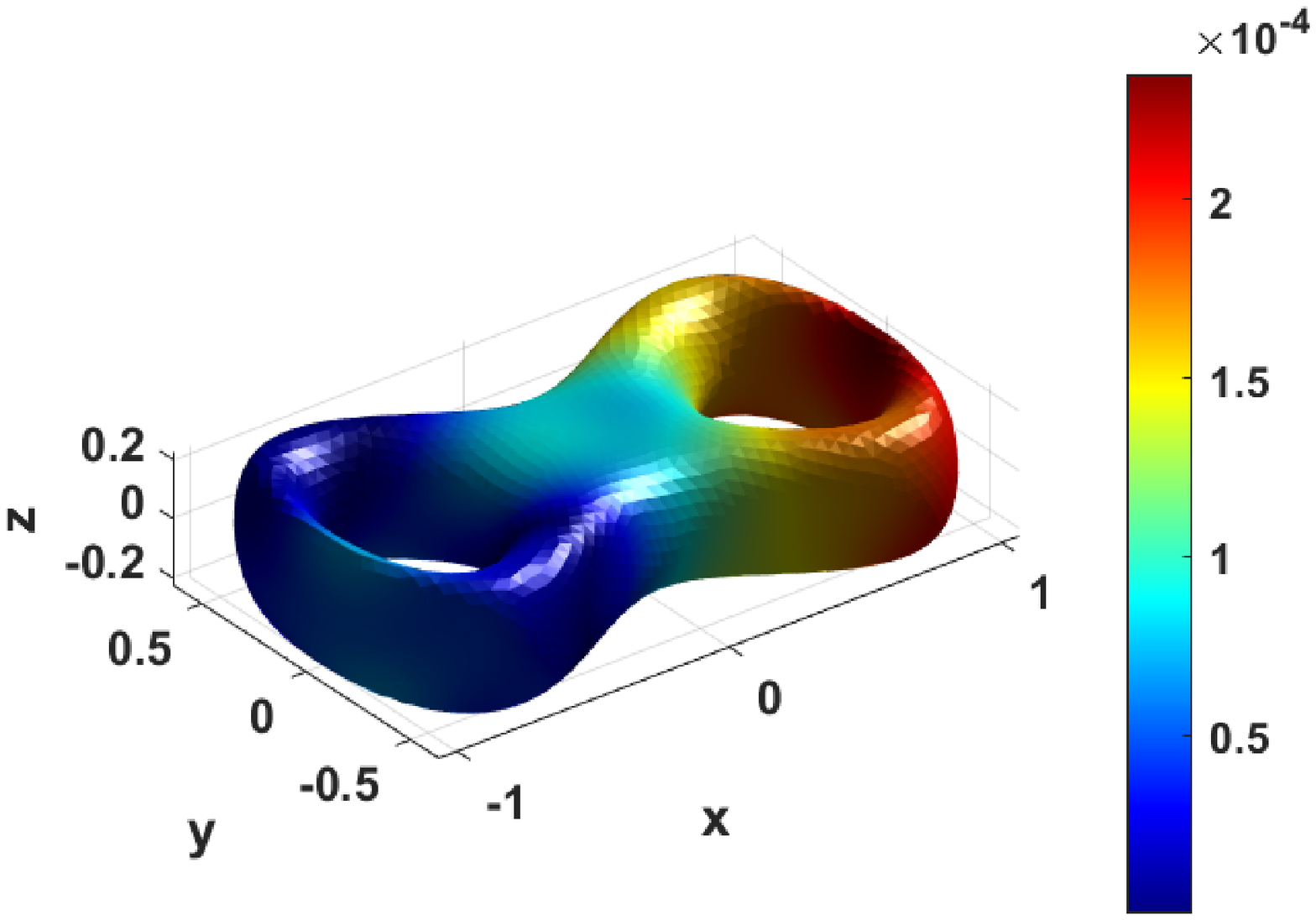}}
\subfigure[]{\label{Fig7c}
\includegraphics[scale=0.4]{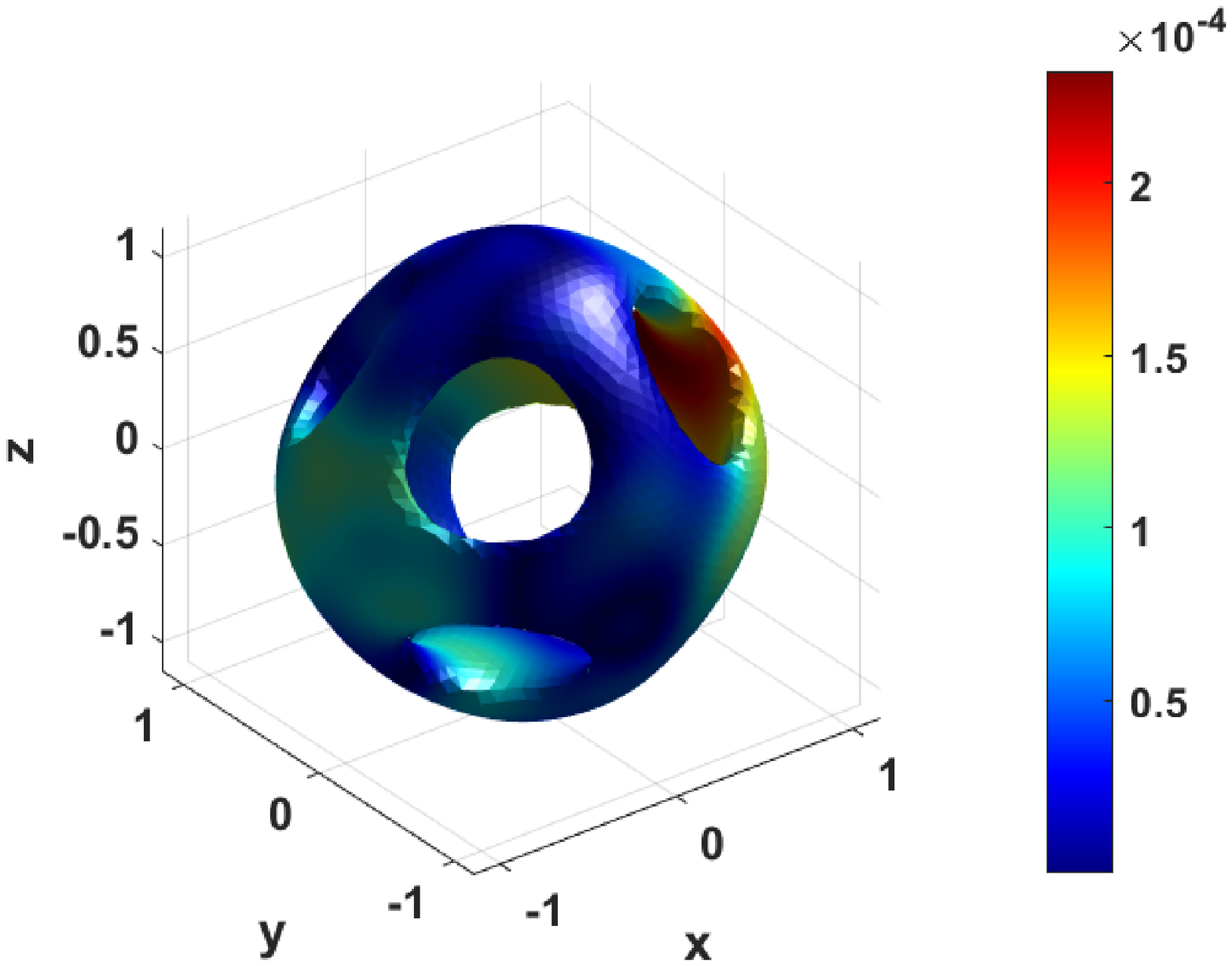}}
\subfigure[]{\label{Fig7d}
\includegraphics[scale=0.4]{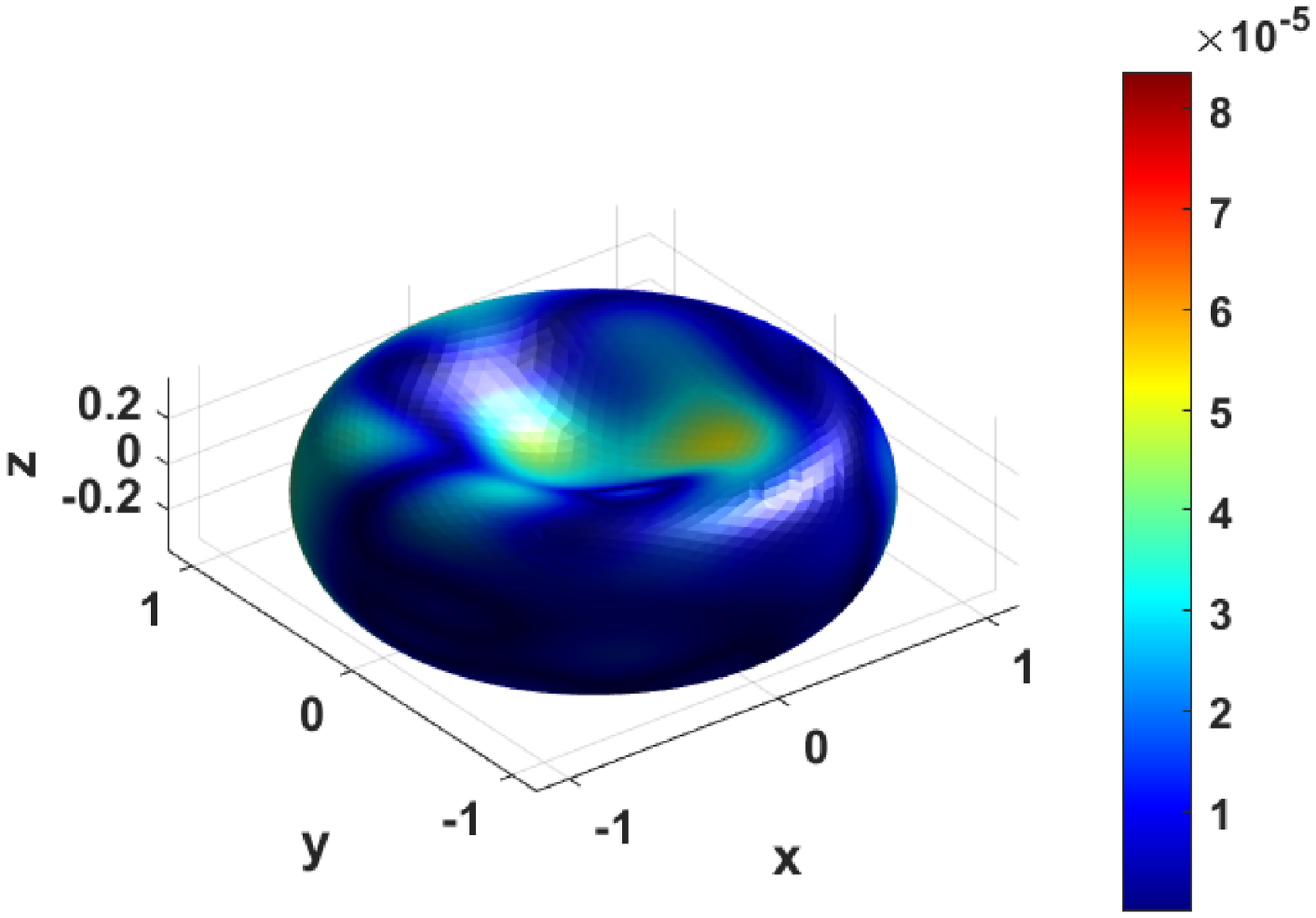}}
\caption{Example 2: Snapshots of error distribution on different manifolds: (a) CDP (b) Bretzel2 (c) Orthocircle (d) RBC.}\label{fig7}
\end{figure}
\begin{table}[!ht]
\begin{center}
\caption{ Example 2: $L_2$ error on different manifolds.}\label{tab1}
\begin{tabular}[]{ccccc}  \hline
Manifolds     & CDP        & Bretzel2   & Orthocircle  & RBC   \\ \hline
$L_2$ error   & 1.18E-03   & 1.51E-03   & 4.20E-03     & 2.37E-03 \\ \hline
\end{tabular}
\end{center}
\end{table}

}
\end{example}

\section{Conclusions}
Physics-informed neural networks (PINNs), as one split-new method belonging to machine learning methods, show a good performance and potential in the solution of the second-order elliptic Partial Differential Equations (PDEs) on 3D manifolds. We could conclude from the first example that PINNs converge rapidly at the beginning of the increasing number of training points due to the dominant effect of the discretization error, while the solution will not be obviously improved with the further increase of the number of training points due to the dominant effect of optimization error. Second example points out that PINNs, as a data-driven method, will not lose accuracy with the dimensionality (shape) increasing in complexity as long as the data provided are accurate enough. This indicates the PINNs have the good stability and robustness. 

This is the first attempt for PINNs in solving PDEs on manifolds. Although PDEs used in this work are relatively simple, it shows a great potential for PINNs in more complicated PDEs (time-dependent equations, nonlinear equations).

\section*{Acknowledgement}
The work described in this paper was supported by the National Science Funds of China (Grant No. 11772119), the Fundamental Research Funds for the Central Universities (Grant No. B200202124), the Foundation for Open Project of State Key Laboratory of Mechanics and Control of Mechanical Structures (Nanjing University Of Aeronautics And Astronautics) (Grant No. MCMS-E-0519G01), the Six Talent Peaks Project in Jiangsu Province of China (Grant No. 2019-KTHY-009) and the Postgraduate Research and Practice Innovation Program of Jiangsu Province (Grant No. KYCX20\_0427).

\section*{Compliance with ethical standards}
\noindent \textbf{Conflict of interests} \ \   On behalf of all authors, the corresponding author states that there is no conflict of interest.

\bibliography{Reference}

\end{document}